\documentclass[12pt]{article}

\usepackage{pstricks,pst-node}
\usepackage{amsmath,amsthm}
\usepackage{amssymb}
\usepackage{color}
\usepackage{xspace}
\usepackage[colorlinks=true,
linkcolor=purple1,
filecolor=brown,
citecolor=blue]{hyperref}

\definecolor{gray1}{rgb}{0.5,0.5,0.5}
\definecolor{black}{rgb}{0,0,0}
\definecolor{purple1}{rgb}{0.2,0,0.4}

\def\gray1{\textcolor{gray1} }

\def\ff{\ensuremath{F}\xspace}

\def\g{\ensuremath{\mathcal G}\xspace}

\def\s{s-partition\xspace}
\def\ss{s-partitions\xspace}
\def\z{x}

\def\red{\textcolor{red} }

\def\v{\vert}
\def\ep{\epsilon}
\def\a{\ensuremath{\mathcal A}\xspace}

\def\fcal{\ensuremath{\mathcal F}\xspace}
\def\q{\ensuremath{\mathcal F}\xspace}
\def\g{\ensuremath{\mathcal G}\xspace}
\def\n{\ensuremath{\mathcal N}\xspace}

\def\p{\ensuremath{\mathcal P}\xspace}

\def\gf{generating function\xspace}

\def\mbf#1{\mathchoice{\hbox{\boldmath $\displaystyle #1$}}
        {\hbox{\boldmath $\textstyle #1$}}
        {\hbox{\boldmath $\scriptstyle #1$}}
        {\hbox{\boldmath $\scriptscriptstyle #1$}}} 

\catcode`\@=11

\thicklines
\newskip\Einheit \Einheit=.6cm
\newcount\xcoord \newcount\ycoord
\newdimen\xdim \newdimen\ydim \newdimen\PfadD@cke \newdimen\Pfadd@cke
\def\PfadDicke#1{\PfadD@cke#1 \divide\PfadD@cke by2 
\Pfadd@cke\PfadD@cke \multiply\PfadD@cke by2}
\long\def\LOOP#1\REPEAT{\def\BODY{#1}\ITERATE}
\def\ITERATE{\BODY \let\next\ITERATE \else\let\next\relax\fi \next}
\let\REPEAT=\fi
\def\Punkt{\hbox{\raise-2pt\hbox to0pt{\hss\scriptsize$\bullet$\hss}}}

\def\DuennPunkt(#1,#2){\unskip
  \raise#2 \Einheit\hbox to0pt{\hskip#1 \Einheit
          \raise-1.5pt\hbox to0pt{\hss\tiny$\bullet$\hss}\hss}}
		  
\def\NormalPunkt(#1,#2){\unskip
  \raise#2 \Einheit\hbox to0pt{\hskip#1 \Einheit
          \raise-3pt\hbox to0pt{\hss\large$\bullet$\hss}\hss}}
\def\DickPunkt(#1,#2){\unskip
  \raise#2 \Einheit\hbox to0pt{\hskip#1 \Einheit
          \raise-4pt\hbox to0pt{\hss\Large$\bullet$\hss}\hss}}
\def\Kreis(#1,#2){\unskip
  \raise#2 \Einheit\hbox to0pt{\hskip#1 \Einheit
          \raise-4pt\hbox to0pt{\hss\Large$\circ$\hss}\hss}}
\def\Diagonale(#1,#2)#3{\unskip\leavevmode
  \xcoord#1\relax \ycoord#2\relax
      \raise\ycoord \Einheit\hbox to0pt{\hskip\xcoord \Einheit
         \unitlength\Einheit
         \line(1,1){#3}\hss}}
\def\AntiDiagonale(#1,#2)#3{\unskip\leavevmode
  \xcoord#1\relax \ycoord#2\relax \advance\xcoord by -0.05\relax
      \raise\ycoord \Einheit\hbox to0pt{\hskip\xcoord \Einheit
         \unitlength\Einheit
         \line(1,-1){#3}\hss}}
\def\Pfad(#1,#2),#3\endPfad{\unskip\leavevmode
  \xcoord#1 \ycoord#2 \thicklines\ZeichnePfad#3\endPfad\thinlines}
\def\ZeichnePfad#1{\ifx#1\endPfad\let\next\relax
  \else\let\next\ZeichnePfad
    \ifnum#1=1
      \raise\ycoord \Einheit\hbox to0pt{\hskip\xcoord \Einheit
         \vrule height\Pfadd@cke width1 \Einheit depth\Pfadd@cke\hss}%
      \advance\xcoord by 1
     \else\ifnum#1=2
      \raise\ycoord \Einheit\hbox to0pt{\hskip\xcoord \Einheit
         \unitlength\Einheit
         \line(0,1){1}\hss}
      \advance\xcoord by 0
      \advance\ycoord by 1
 \else\ifnum#1=3
      \raise\ycoord \Einheit\hbox to0pt{\hskip\xcoord \Einheit
         \unitlength\Einheit
         \line(1,1){1}\hss}
      \advance\xcoord by 1
      \advance\ycoord by 1
    \else\ifnum#1=4
      \raise\ycoord \Einheit\hbox to0pt{\hskip\xcoord \Einheit
         \unitlength\Einheit
         \line(1,-1){1}\hss}
      \advance\xcoord by 1
      \advance\ycoord by -1
   \else\ifnum#1=5
      \raise\ycoord \Einheit\hbox to0pt{\hskip\xcoord \Einheit
         \unitlength\Einheit
         \line(2,1){2}\hss}
      \advance\xcoord by 2
      \advance\ycoord by 1
	  \else\ifnum#1=6
      \raise\ycoord \Einheit\hbox to0pt{\hskip\xcoord \Einheit
         \unitlength\Einheit
         \line(2,-1){2}\hss}
      \advance\xcoord by 2
      \advance\ycoord by -1
	  \else\ifnum#1=7
      \raise\ycoord \Einheit\hbox to0pt{\hskip\xcoord \Einheit
         \unitlength\Einheit
         \line(3,1){3}\hss}
      \advance\xcoord by 3
      \advance\ycoord by 1
	  \else\ifnum#1=8
      \raise\ycoord \Einheit\hbox to0pt{\hskip\xcoord \Einheit
         \unitlength\Einheit
         \line(3,-1){3}\hss}
      \advance\xcoord by 3
      \advance\ycoord by -1
    \fi\fi\fi\fi\fi\fi\fi\fi
  \fi\next}
\def\hSSchritt{\leavevmode\raise-.4pt\hbox 
to0pt{\hss.\hss}\hskip.2\Einheit
  \raise-.4pt\hbox to0pt{\hss.\hss}\hskip.2\Einheit
  \raise-.4pt\hbox to0pt{\hss.\hss}\hskip.2\Einheit
  \raise-.4pt\hbox to0pt{\hss.\hss}\hskip.2\Einheit
  \raise-.4pt\hbox to0pt{\hss.\hss}\hskip.2\Einheit}
\def\vSSchritt{\vbox{\baselineskip.2\Einheit\lineskiplimit0pt
\hbox{.}\hbox{.}\hbox{.}\hbox{.}\hbox{.}}}
\def\DSSchritt{\leavevmode\raise-.4pt\hbox to0pt{%
  \hbox to0pt{\hss.\hss}\hskip.2\Einheit
  \raise.2\Einheit\hbox to0pt{\hss.\hss}\hskip.2\Einheit
  \raise.4\Einheit\hbox to0pt{\hss.\hss}\hskip.2\Einheit
  \raise.6\Einheit\hbox to0pt{\hss.\hss}\hskip.2\Einheit
  \raise.8\Einheit\hbox to0pt{\hss.\hss}\hss}}
\def\dSSchritt{\leavevmode\raise-.4pt\hbox to0pt{%
  \hbox to0pt{\hss.\hss}\hskip.2\Einheit
  \raise-.2\Einheit\hbox to0pt{\hss.\hss}\hskip.2\Einheit
  \raise-.4\Einheit\hbox to0pt{\hss.\hss}\hskip.2\Einheit
  \raise-.6\Einheit\hbox to0pt{\hss.\hss}\hskip.2\Einheit
  \raise-.8\Einheit\hbox to0pt{\hss.\hss}\hss}}
\def\SPfad(#1,#2),#3\endSPfad{\unskip\leavevmode
  \xcoord#1 \ycoord#2 \ZeichneSPfad#3\endSPfad}
\def\ZeichneSPfad#1{\ifx#1\endSPfad\let\next\relax
  \else\let\next\ZeichneSPfad
    \ifnum#1=1
      \raise\ycoord \Einheit\hbox to0pt{\hskip\xcoord \Einheit
         \hSSchritt\hss}%
      \advance\xcoord by 1
    \else\ifnum#1=2
      \raise\ycoord \Einheit\hbox to0pt{\hskip\xcoord \Einheit
        \hbox{\hskip-2pt \vSSchritt}\hss}%
      \advance\ycoord by 1
    \else\ifnum#1=3
      \raise\ycoord \Einheit\hbox to0pt{\hskip\xcoord \Einheit
         \DSSchritt\hss}
      \advance\xcoord by 1
      \advance\ycoord by 1
    \else\ifnum#1=4
      \raise\ycoord \Einheit\hbox to0pt{\hskip\xcoord \Einheit
         \dSSchritt\hss}
      \advance\xcoord by 1
      \advance\ycoord by -1
    \fi\fi\fi\fi
  \fi\next}
\def\Koordinatenachsen(#1,#2){\unskip
 \hbox to0pt{\hskip-.5pt\vrule height#2 \Einheit width.5pt depth1 
\Einheit}%
 \hbox to0pt{\hskip-1 \Einheit \xcoord#1 \advance\xcoord by1
    \vrule height0.25pt width\xcoord \Einheit depth0.25pt\hss}}
\def\Koordinatenachsen(#1,#2)(#3,#4){\unskip
 \hbox to0pt{\hskip-.5pt \ycoord-#4 \advance\ycoord by1
    \vrule height#2 \Einheit width.5pt depth\ycoord \Einheit}%
 \hbox to0pt{\hskip-1 \Einheit \hskip#3\Einheit 
    \xcoord#1 \advance\xcoord by1 \advance\xcoord by-#3 
    \vrule height0.25pt width\xcoord \Einheit depth0.25pt\hss}}
\def\Gitter(#1,#2){\unskip \xcoord0 \ycoord0 \leavevmode
  \LOOP\ifnum\ycoord<#2
    \loop\ifnum\xcoord<#1
      \raise\ycoord \Einheit\hbox to0pt{\hskip\xcoord 
\Einheit\Punkt\hss}%
      \advance\xcoord by1
    \repeat
    \xcoord0
    \advance\ycoord by1
  \REPEAT}
\def\Gitter(#1,#2)(#3,#4){\unskip \xcoord#3 \ycoord#4 \leavevmode
  \LOOP\ifnum\ycoord<#2
    \loop\ifnum\xcoord<#1
      \raise\ycoord \Einheit\hbox to0pt{\hskip\xcoord 
\Einheit\Punkt\hss}%
      \advance\xcoord by1
    \repeat
    \xcoord#3
    \advance\ycoord by1
  \REPEAT}
\def\Label#1#2(#3,#4){\unskip \xdim#3 \Einheit \ydim#4 \Einheit
  \def\lo{\advance\xdim by-.5 \Einheit \advance\ydim by.5 \Einheit}%
  \def\llo{\advance\xdim by-.25cm \advance\ydim by.5 \Einheit}%
  \def\loo{\advance\xdim by-.5 \Einheit \advance\ydim by.25cm}%
  \def\o{\advance\ydim by.25cm}%
  \def\ro{\advance\xdim by.5 \Einheit \advance\ydim by.5 \Einheit}%
  \def\rro{\advance\xdim by.25cm \advance\ydim by.5 \Einheit}%
  \def\roo{\advance\xdim by.5 \Einheit \advance\ydim by.25cm}%
  \def\l{\advance\xdim by-.30cm}%
  \def\r{\advance\xdim by.30cm}%
  \def\lu{\advance\xdim by-.5 \Einheit \advance\ydim by-.6 \Einheit}%
  \def\llu{\advance\xdim by-.25cm \advance\ydim by-.6 \Einheit}%
  \def\luu{\advance\xdim by-.5 \Einheit \advance\ydim by-.30cm}%
  \def\u{\advance\ydim by-.30cm}%
  \def\ru{\advance\xdim by.5 \Einheit \advance\ydim by-.6 \Einheit}%
  \def\rru{\advance\xdim by.25cm \advance\ydim by-.6 \Einheit}%
  \def\ruu{\advance\xdim by.5 \Einheit \advance\ydim by-.30cm}%
  #1\raise\ydim\hbox to0pt{\hskip\xdim
     \vbox to0pt{\vss\hbox to0pt{\hss$#2$\hss}\vss}\hss}%
}
\catcode`\@=12

\begin{document}
\newtheorem{theorem}{Theorem}
\newtheorem{defn}[theorem]{Definition}
\newtheorem{lemma}[theorem]{Lemma}
\newtheorem{prop}[theorem]{Proposition}
\newtheorem{cor}[theorem]{Corollary}
\newtheorem{conj}[theorem]{Conjecture}
\newtheorem*{mainthm}{Main Theorem}
\newtheorem*{gnrl}{Main Theorem Generalized}

\begin{center}
{\Large
The Run Transform   \\ 
}

\vspace{10mm}

\begingroup
\renewcommand\thefootnote{\fcalnsymbol{footnote}}
\begin{tabular}{l@{\hspace{3em}}l}
   David Callan  &   Emeric Deutsch \\[1mm]
  Department of Statistics &  Department of Mathematics\\[-1mm]
  University of Wisconsin-Madison &  Polytechnic Institute of NYU\\[-1mm]
   Madison, WI \ 53706 &   Brooklyn, NY 11201 \\[1mm]
  {\small {\tt callan@stat.wisc.edu}}   & {\small {\tt emericdeutsch@msn.com}} 
\end{tabular}
\endgroup

\end{center}

\begin{abstract}
We consider the transform from sequences to triangular arrays defined 
in terms of generating functions by $f(x) \rightarrow 
\frac{1-x}{1-xy} f\big(\frac{x(1-x)}{1-xy}\big)$. We establish a 
criterion for the transform of a nonnegative sequence to be 
nonnegative, and we show that the transform counts certain classes of lattice paths 
by number of ``pyramid ascents'', as well as certain classes of ordered partitions by number of blocks that consist of increasing consecutive integers.
\end{abstract}

\section{Introduction} 
\vspace*{-2mm}
We investigate the transform $\Phi$ defined on formal power series $f(x)$ 
by
\[
\Phi\big(f(x)\big) :=
\frac{1-x}{1-xy}\: f\!\left(\frac{x\,(1-x)}{1-xy}\right).
\]
Following Herbert Wilf's dictum, 
``A generating function is a clothesline on which we hang up a sequence
of numbers for display'' \cite[Chapter 1]{wilf}, we will use sequences/arrays 
and their generating functions interchangeably. Thus the transform
$\Phi$ is also defined for sequences $(a_{n})_{n\ge 0}$. It turns out 
that the transform is closely related to the Catalan numbers and 
there is a nice combinatorial interpretation for the transform of the 
size-counting sequence for various classes of partitions into sets of 
lists (blocks) and  various classes of lattice paths of upsteps $U$, 
flatsteps $F$, and downsteps $D$. In the former case, the transform counts partitions by 
number of runs, where a \emph{run}, also known as an \emph{adjacent block} \cite{adj}, 
is a block that 
consists of increasing consecutive integers. In the latter case it 
counts lattice paths by number of pyramid ascents, where an \emph{ascent} is a 
maximal subpath of the form $U^{k},\ k\ge 1$, a \emph{pyramid} is a 
maximal subpath of the form $U^{k}D^{k}, k\ge 1$, and a \emph{pyramid ascent} 
is an ascent that is the first half of a pyramid. For example, among the four 
ascents of $UUDUUDUUDDDUDD$, only the last two
($UU$ and $U$) are pyramid ascents. 

Because of the interpretation in terms of runs, and for brevity, we will designate 
$\Phi$ the \emph{run transform}.

In Section \ref{catconv} we review the Catalan numbers and two of 
their interpretations, and in Section 3 we establish basic properties of the run transform. 
Section 4 gives a criterion for the run transform 
of a nonnegative sequence to also be nonnegative. 
Section 5 gives interpretations of the run transform of the Catalan 
numbers in terms of both Dyck paths and noncrossing partitions, the 
basis for subsequent generalizations.
Section 6 generalizes to paths of $j$-upsteps $(j,j)$ and downsteps 
$(1,-1)$.
Section 7 recalls the notion of a set-of-lists partition, \s for 
short, and introduces the notion of a run-closed family of \ss  and states 
the result 
that if $f(\z)$ is the \gf by size of a run-closed family $\fcal$ of \ss, 
then the  run transform of $f(\z)$ counts $\fcal$ by size and number of runs.
This result is proved in 
Sections 8 and 9, and generalized in Section 10.
Section 11 considers the transform for paths of upsteps, flatsteps, 
and downsteps and Section 12 presents a conjecture.

Sequences in The On-Line Encyclopedia of Integer Sequences (OEIS) \cite{oeis}, 
are referred to by their six-digit A-numbers.

\section{Review of the Catalan numbers, Dyck paths, and noncrossing 
partitions}\vspace*{-5mm}\label{catconv}
The Catalan numbers 
(sequence \htmladdnormallink{A000108}{http://oeis.org/A000108} in OEIS)
are intimately related to the run transform $\Phi$; so let us recall 
some facts and fix some notation for them and for two of their 
interpretations. The \gf for the Catalan 
numbers $C_{n}=\frac{1}{n+1}\binom{2n}{n}=\binom{2n}{n}-\binom{2n}{n-1}$ is 
$C(x)=(1-\sqrt{1-4x})/(2x)$. 
The Catalan convolution matrix is defined by 
$C=\left(\binom{2j-i}{j-i} - \binom{2j-i}{j-i-1} \right)_{i,j\ge 0}$ and its 
inverse is given by $C^{-1}=
\left((-1)^{j-i}\binom{i+1}{j-i}\right)_{i,j\ge 0}$
The first few rows and columns are shown.
\[
C=
\left(\begin{array}{cccccccc}
    1 & 1 & 2 & 5 & 14 & 42 & 132 & \ldots  \\
     & 1 & 2 & 5 & 14 & 42 & 132 & \ldots \\
     &  & 1 & 3 & 9 & 28 & 90 & \ldots \\
     &  &  & 1 & 4 & 14 & 48 & \ldots \\
     &  &  &  & 1 & 5 & 20 & \ldots \\
     &  &  &  &  & 1 & 6 & \ldots  \\
     &  &  &  &  &  & 1 & \ldots \\
      &  &  &  &  &  &   & \ddots 
\end{array}\right),\quad C^{-1}=
\left(\begin{array}{cccccccc}
    1 & -1 &  0 & 0  & 0  & 0  & 0 & \ldots  \\
     & 1 & -2 & 1 & 0  & 0  & 0 & \ldots \\
     &  & 1 & -3 & 3 & -1 & 0  & \ldots \\
     &  &  & 1 & -4 & 6 & -4 & \ldots \\
     &  &  &  & 1 & -5 & 10 & \ldots \\
     &  &  &  &  & 1 & -6 & \ldots  \\
     &  &  &  &  &  & 1 & \ldots \\
     &  &  &  &  &  &   & \ddots 
\end{array}\right)
\]
Thus the top row (row 0) of $C$, the Catalan numbers, has \gf $C(x)$.
It is well known that the entries of $C$  count nonnegative lattice 
paths. Specifically, define a $U$-$D$ path to be a lattice path of \emph{upsteps} $U=(1,1)$ and 
\emph{downsteps} $D=(1,-1)$ (not necessarily starting at the origin). Let $\n_{mn}\ (=\n_{m,n})$ denote the set of $U$-$D$
paths consisting of $n$ upsteps and $m+n$ downsteps that never dip below ground 
level, the horizontal line through the terminal vertex. The size of a $U$-$D$
path is its number of downsteps. 
Then $\v\,\n_{mn}\,\v = C_{mn}$.
Set  $\n_{m}:=\bigcup_{n\ge 0}\n_{mn}$. A \emph{Dyck} path is a member 
of $\n_{0}$.
It is also well known 
(see, e.g., \cite[Remark 5]{woan} or \cite{tedford}) that row $m$ of $C$, starting at 
the diagonal entry, is the $(m+1)$-fold convolution of the top row. 
Hence, the \gf for $\n_{m}$ by size is $x^{m}C(x)^{m+1}$, and,
with $\mathbf{x}$ defined to be the column vector 
$(1,x,x^{2},\ldots)^{\,t}$, we have the matrix-vector product
\begin{equation}\label{matvec}
C\,\mathbf{x}=(C(x),\:xC(x)^{2},\:x^{2}C(x)^{3},\:\ldots)^{\,t}.
\end{equation} 
We define $C(x,y)$ to be $\Phi\big(C(x)\big)$. Thus
\begin{equation}\label{ctz}
C(x,y)=\frac{1-\sqrt{1-4\,\frac{\textrm{{\small $x\,(1-x)$}}}{\textrm{{\small 
$1-xy$}}}}}{2x}.
\end{equation} 

The \emph{matching step} of a given step in a Dyck 
path is the other end-step of the shortest Dyck subpath containing 
the given step as an end-step.
Ascents (and pyramid ascents) were defined in the Introduction and, analogously, 
a \emph{descent} is a maximal subpath of the form $D^{k},\ k\ge 1$.

A nonempty Dyck path decomposes (at its 
returns to ground level) into \emph{components}, each of which is a \emph{primitive} 
Dyck path---a nonempty Dyck path 
whose only return to ground level is at the end.

There is a well known bijection from Dyck paths to noncrossing 
partitions, due to Simion \cite{nc2000}.
Traverse the Dyck path from left to right and number the down steps from 
1 to $n$.
Give the same labels to the matching up steps. The numbers
on the ascents form the blocks of the partition.
Under this bijection, pyramid ascents become runs in the partition.

\section{Basic properties of the run transform}\vspace*{-5mm}
We will use $F(x,y)$ for $\Phi\big(f(x)\big)$ to show the dependency on both variables.
Clearly, $F(x,1)=f(x)$ and so the row sums of the run transform give the 
original sequence. The run transform $\Phi$ is linear, 
\[
\Phi\big(\alpha f(x)+\beta g(x)\big)=
\alpha \Phi\big(f(x)\big)+\beta \Phi\big(g(x)\big),
\]
and has a multiplicativity property, \[
\Phi\big(xf(x)g(x)\big)=x\Phi\big(f(x)\big)\Phi\big(g(x)\big).
\]
 More generally,
\[
\Phi\big(x^{i-1}f_{1}(x)f_{2}(x)... f_{i}(x)\big)=x^{i-1}\Phi\big(f_{1}(x)\big)
\Phi\big(f_{2}(x)\big)... \Phi\big(f_{i}(x)\big).
\]
 In particular, for $i=k+1$ and $f=f_{1}=f_{2}=...,$
\[
  \Phi\big(x^{k}f(x)^{k+1}\big)=x^{k}\Phi\big(f(x)\big)^{k+1}.
  \]
From this fact, together with linearity, we obtain
\begin{lemma} \label{prop1}
    For an arbitrary sequence $(a_{k})_{k\ge 0}$, 
\[
\Phi\left(\sum_{k\ge 0}a_{k}x^{k}f(x)^{k+1}\right)=
     \sum_{k\ge 0}a_{k}x^{k}\Phi\big(f(x)\big)^{k+1}.
\]
\end{lemma} 
\qed

\begin{prop} \label{em}
    Let $\mathbf{b}=(b_{k})_{k\ge 0}$ be an arbitrary sequence. 
    Then its  run transform  is
    \[
    \sum_{k\ge 0}a_{k}x^{k}C(x,y)^{k+1},
    \]
    where $\mathbf{a}=(a_{k})_{k\ge 0}$ is defined by 
    $\mathbf{a}=\mathbf{b}\,C^{-1}$.
\end{prop}    

Proof. The defining relation  $\mathbf{a}=\mathbf{b}\,C^{-1}$ yields 
 $\mathbf{b}=\mathbf{a}\,C$ and, multiplying by the column vector  
 $\mathbf{x}$, 
 \[
 \mathbf{b}\,\mathbf{x}=\mathbf{a}\,C\,\mathbf{x}
 \]
 which, making use of (\ref{matvec}), translates into 
 \[
 f(x)=\sum_{k\ge 0}a_{k}x^{k}C(x)^{k+1}.
 \]
Now apply Lemma \ref{prop1} with $f(x)=C(x)$. \qed

\section{A criterion for nonnegativity}
\vspace*{-5mm}
\begin{prop}
For a nonnegative sequence $\mathbf{a}=(a_{k})_{k\ge 0}$, its 
run transform is nonnegative if and only if the sequence 
$\mathbf{x}:=\mathbf{a}\,C^{-1}$ is nonnegative. 
\end{prop}    
Proof. We have $C(x,0)=(1-x)C\big(x(1-x)\big)=1$ and it follows from  
Proposition \ref{em} that column 0 of the run transform, given by $F(x,0)$, is (the transpose of)  
$\mathbf{x}$. So the condition is certainly necessary. 
Sufficiency will follow if we know that each power of $C(x,y)$ is the 
\gf of a nonnegative array. 
For $C(x,y)$ itself, nonnegativity follows from a combinatorial interpretation in terms of 
decorated forests \cite[Section 9]{lcoforest} or from Proposition \ref{pyrascDyck} below, 
but we can also give an analytic proof as follows. Say $(u_{i,j})_{i\ge 0,\ 0\le j \le i}$ is the array 
of coefficients for $C(x,y)$. We have the identity $(2xC(x,y)-1)^{2}=1-4x(1-x)/(1-xy)$, 
leading to
\[
xC(x,y)^{2}=C(x,y)-\frac{1-x}{1-xy}.
\]
Picking out coefficients leads to a recurrence for $u_{i,j}:\ u_{0,0} = 1$, and 
\[
u_{n,k}  = 
\begin{cases}
     \sum_{i=0}^{n-1}\sum_{j=0}^{n}u_{i,j}u_{n-1-i,n-j} +1, & \textrm{  if $1\le k \le n$;} \\
     \sum_{i=0}^{n-1}\sum_{j=0}^{n-1}u_{i,j}u_{n-1-i,n-1-j} -1  & \textrm{  if $0\le k = n-1$;} \\
     \sum_{i=0}^{n-1}\sum_{j=0}^{k}u_{i,j}u_{n-1-i,k-j}  & \textrm{  if $0\le k \le n-2$,}
\end{cases}  
\]
from which it is easy to see that $u_{i,j}$ is nonnegative, the $-1$ in the middle equality 
notwithstanding. Finally, it is 
easy to check that nonnegativity of $C(x,y)$ implies nonnegativity 
of all its powers.  

\section{The run transform of the Catalan numbers}
\begin{lemma}\label{pyrascDyck}
    The run transform $C(x,y)$ of the Catalan number \gf $C(x)$ counts 
    Dyck paths by size and number of pyramid ascents, equivalently, 
    noncrossing partitions by size and number of runs.
\end{lemma}
Proof. Let $F(x,y)$ denote the \gf for Dyck paths by size and number of pyramid ascents. 
A Dyck path $P$ is either empty or has the decomposition $P=U^{r}DP_{1}DP_{2}D... DP_{r}$ 
for some $r\ge1$, where the $P_{i}$ are Dyck paths. 
Each pyramid ascent in $P_{1},...,P_{r}$ is a pyramid ascent in $P$ and, 
if $P_{1},...,P_{r-1}$ are all empty paths, then the first ascent of $P$ 
is also a pyramid ascent, contributing a $y$ factor. We thus obtain
\[
F=1+\sum_{r\ge 1}x^{r}\big(y+F^{r-1}-1\big)F
\]
which leads at once to
\[
xF^{2}-F+\frac{1-x}{1-xy} = 0,
\]
an equation whose solution is $F(x,y)=C(x,y)$. \qed

Recall that  $\n_{k}$ is the set of nonnegative $U$-$D$ paths with 
$k$ more  downsteps than upsteps.
\begin{theorem}
    The run transform  of the \gf for $\n_{k}$ by size counts 
    $\n_{k}$ by size and number of pyramid ascents.
\end{theorem}
Proof. A path in $\n_{k}$ decomposes as $P_{1}DP_{2}D\ldots 
P_{k}DP_{k+1}$ with each $P_{i}$ a Dyck path. The statistics size and 
number of pyramid ascents are additive over this decomposition. So one 
multiplies the generating functions given by Proposition 
\ref{pyrascDyck}, and 
the \gf for $\n_{k}$ by size and number of pyramid 
ascents is $x^{k}C(x,y)^{k+1}$. By Lemma 1, the run transform of 
$x^{k}C(x)^{k+1}$ is $x^{k}C(x,y)^{k+1}$. \qed

In subsequent sections we generalize this result in 3 ways: (i) from Dyck 
paths to 
$U$-$D$ paths in which each ascent has length divisible by $j$, 
(ii) from Dyck 
paths to $U-F-D$ paths, where flatsteps $F=(1,0)$ are allowed, (iii) 
from  noncrossing partitions to run-closed families of \ss, 
defined in Section \ref{runclosed} below.

\section[U\^{}\{j\}-D paths]{$\mbf{U^{j}}$-$\mbf{D}$ paths}
\vspace*{-5mm}
Fix a positive integer $j$. 
A \emph{$j$-Dyck} path is a Dyck path in which each ascent has length 
divisible by $j$. Equivalently, it can be viewed as a nonnegative 
lattice path of so-called $j$-upsteps $(j,j)$ and (ordinary) downsteps 
$(1,-1)$. Its \emph{size} is the number of $j$-upsteps, equivalently, 
(number of downsteps)$/j$. A \emph{$j$-nice} pyramid ascent 
is one that ends at height $\equiv 0$ (mod $j$).

\begin{lemma}\label{jud}
    The run transform  of the \gf for $j$-Dyck paths by size is the \gf for $j$-Dyck paths 
    by size and number of pyramid ascents.
\end{lemma}
Proof. 
Let $\ff_{0}(x,y)$ denote the \gf for $j$-Dyck 
paths with $x$ marking size and $y$ marking the number of $j$-nice pyramid 
ascents. To find an equation for $\ff_{0}(x,y)$, it is convenient to introduce
the \gf $\ff_{i}(x,y)$, $i$ an integer, for the number of pyramid 
ascents that end at height $\equiv i$ (mod $j$). 
Set $\ff = \ff_{0}\ff_{1}...\ff_{j-1}$.

A nonempty $j$-Dyck path $P$ has a decomposition as illustrated for $j=3$ and $r=2$ where 
$r$ is the number of initial $j$-upsteps, and the $P_{i}$'s are all $j$-Dyck paths (possibly empty).
\Einheit=0.5cm
\[
\Pfad(-12,0),3333334\endPfad
\Pfad(-3,5),4\endPfad
\Pfad(0,4),4\endPfad
\Pfad(3,3),4\endPfad
\Pfad(6,2),4\endPfad
\Pfad(9,1),4\endPfad
\SPfad(-12,0),111111111111111111111111\endSPfad
\Label\o{P_{1}}(-4,5)
\Label\o{P_{2}}(-1,4)
\Label\o{P_{3}}(2,3)
\Label\o{P_{4}}(5,2)
\Label\o{P_{5}}(8,1)
\Label\o{P_{6}}(11,0)
\DuennPunkt(-11,1)
\DuennPunkt(-10,2)
\DuennPunkt(-8,4)
\DuennPunkt(-7,5)
\NormalPunkt(-12,0)
\NormalPunkt(-9,3)
\NormalPunkt(-6,6)
\NormalPunkt(-5,5)
\NormalPunkt(-3,5)
\NormalPunkt(-2,4)
\NormalPunkt(0,4)
\NormalPunkt(1,3)
\NormalPunkt(3,3)
\NormalPunkt(4,2)
\NormalPunkt(6,2)
\NormalPunkt(7,1)
\NormalPunkt(9,1)
\NormalPunkt(10,0)
\NormalPunkt(12,0)
\Label\u{\textrm{A $j$-Dyck path with $j=3$}}(0,-1)
\]

\noindent In general, the decomposition is
\[
U^{jr}DP_{1}DP_{2}... DP_{jr},
\]
for some $r\ge 1$.
To obtain an expression for $\ff_{0}$ from this decomposition, $P_{1}$ contributes $\ff_{1}$ 
because it starts at height $\equiv -1$ mod $j$, 
$P_{2}$ contributes $\ff_{2}$,\ ..., $P_{j}$ contributes $\ff_{0}$, $P_{j+1}$ 
contributes $\ff_{1}$, and so on, cyclically. And if $P_{1},P_{2},...P_{jr-1}$ are all empty, 
the first ascent of $P$ is a pyramid ascent. Splitting into the two cases where 
$P_{1},P_{2},...,P_{jr-1}$ are all empty or not, we thus obtain
\[
\ff_{0}=1+\sum_{r\ge 1}x^{r}\big(y+\ff_{1}\ff_{2}... \ff_{j-1}\ff^{r-1}-1\big)\ff_{0}
\]
which simplifies to
\begin{equation}\label{f0}
\ff_{0}=1+(y-1)\ff_{0} \frac{x}{1-x} +\frac{x\ff}{1-x\ff}.
\end{equation}
For $i\not\equiv 0$ (mod $j$), there is no need to split into cases. We find 
$\ff_{i}=1+x\ff + x^{2}\ff^{2}+...$, leading to 
\begin{equation}\label{fi}
\hspace*{50mm} \ff_{i}=\frac{1}{1-x\ff} \hspace*{30mm}\textrm{ for $i\not\equiv 0$ (mod $j$).}
\end{equation}
Eliminating $f$ from (\ref{f0}) and (\ref{fi}), we obtain
\begin{equation}
  \hspace*{40mm}  \ff_{i}(x,y)=\frac{1-xy}{1-x}\ff_{0}(x,y)\hspace*{20mm}\textrm{ for $i\not\equiv 0$ (mod $j$).}
    \label{fij}
\end{equation}
Hence,
\[
\ff=\left(\frac{1-xy}{1-x}\right)^{j-1}\ff_{0}^{j}
\]
and (\ref{f0}) becomes, after further simplification
\begin{equation}
    (1-x)^{j}-(1-x)^{j-1}(1-xy)\ff_{0} + x(1-xy)^{j}\ff_{0}^{j+1}=0.
    \label{eq:definej}
\end{equation}
From (\ref{eq:definej}), $f(x):=\ff_{0}(x,1)$ has defining 
equation $1-f+xf^{j+1}=0$, and the run transform of $f(x)$ satisfies 
(\ref{eq:definej}). Hence the run transform of $f(x)$ is $\ff_{0}(x,y)$. 
\qed

Now fix nonnegative integers $m$ and $d$.  
A $(j,m,d)$-$U$-$D$ path is a path of $j$-upsteps and downsteps  
that starts (for convenience) at $(0,m)$, has lowest point at level $-d$, and ends on the $x$-axis.
A $(j,0,0)$-$U$-$D$ path is just a $j$-Dyck path.
The size of a $(j,m,d)$-$U$-$D$ path is $\lfloor \textrm{(number of downsteps)}/j \rfloor$, 
and so it is convenient to express $m$ as $jk+\ell$ with $0\le \ell \le j-1$.

\begin{theorem}\label{jmdUD}
    The run transform  of the \gf for $(j,m,d)$-$U$-$D$ paths by size is the \gf for $(j,m,d)$-$U$-$D$ paths 
    by size and number of pyramid ascents.
\end{theorem}
Proof. 
Let $G(x,y)$ denote the \gf for $(j,m,d)$-$U$-$D$ paths with $x$ marking 
size and $y$ marking the number of $j$-nice pyramid ascents. 
A $(j,m,d)$-$U$-$D$ path  has the decomposition:
\begin{equation}
    P_{0}\,D\,P_{1}\,D\,P_{2}\,\ldots\, 
    D\,P_{m+d}\,U^{j}\,P_{m+d-1}'\,U^{j}\,P_{m+d-2}'\, \ldots \,U^{j}\,P_{m}'.
    \label{decomp}
\end{equation}
where the $P_{i}$ and the $P_{i}'$ are all $j$-Dyck paths.

If a $P_{i}'$ begins with a pyramid, the pyramid is killed by the 
immediately preceding $U^{j}$. This necessitates introducing the \gf $H(x,y)$ for $j$-Dyck 
paths that begin with a pyramid:
\[
H(x,y)=\sum_{i\ge 1}x^{i}y\ff_{0} =\frac{xy}{1-x}\ff_{0}.
\]
The decomposition (\ref{decomp}) together with (\ref{fij}) yields
\begin{eqnarray}
  G(x,y)   & = &  x^{k+d}\ff_{-m}\,\ldots\,\ff_{-1}\,\ff_{0}\,
  \ff_{1}\,\ff_{2}\,\ldots\,\ff_{jd}\left(\frac{H}{y}+\ff_{0}-H\right)  \\
     & = & x^{k+d}\ff_{0}^{k+d+1}\ff_{1}^{m-k+(j-1)d}
     \left(\frac{1-xy}{1-x}\ff_{0}\right)^{d}  \\
     & = & x^{k+d}
     \left(\frac{1-xy}{1-x}\right)^{m-k+jd}\ff_{0}^{m+1+(j+1)d}.
\end{eqnarray}

Hence $g(x):=G(x,1)=x^{k+d}\ff_{0}(x,1)^{m+1+(j+1)d}$ is the \gf for $(j,m,d)$-$U$-$D$ paths by size.
 The run transform of $g(x)$ is 
\begin{eqnarray*}
    \frac{1-x}{1-xy}\,g\left(\frac{x(1-x)}{1-xy}\right) & = & 
    \frac{1-x}{1-xy} \left(\frac{x(1-x)}{1-xy}\right)^{k+d}
    \ff_{0}\Big(\frac{x(1-x)}{1-xy},1\Big)^{m+1+(j+1)d} \\
     & = &  x^{k+d}\left(\frac{1-xy}{1-x}\right)^{m+jd-k}
     \left(\frac{1-x}{1-xy}\, \ff_{0}\Big(\frac{x(1-x)}{1-xy},1\Big)\right)^{m+1+(j+1)d}     \\
     & = & x^{k+d}\left(\frac{1-xy}{1-x}\right)^{m+jd-k} 
     \ff_{0}(x,y)^{m+1+(j+1)d} \\
      & = & G(x,y),
\end{eqnarray*}
the next to last equality using the fact that the run transform of $\ff_{0}(x,1)$ 
is $\ff_{0}(x,y)$ (Lemma \ref{jud}).

\section{Run-closed families of \ss}\label{runclosed}
\vspace*{-5mm}
A \emph{set-of-lists} partition, or \emph{\s} for short, also known as a \emph{fragmented 
permutation} \cite[p.\,125]{flaj09}, is a partition $\pi$ of a set $S$ into a set of lists. 
The size of $\pi$, denoted $\v\, \pi \,\v$, is $\v\, S \, \v$.
An \s is \emph{standard} if its support set is an initial segment of the positive integers. 
We use the familiar term \emph{blocks} for the lists in an \s, and we always arrange 
the blocks in increasing order of their first entry. 
Recall that a run is a block that consists  
of consecutive integers in increasing order. Thus the \s 
$3\,8\,1\,/\,4\,5\,6\,/\,7\,2\,/\, 9$ has size 9 and four blocks, two of which 
are runs, $4\,5\,6$ and 9. 
A permutation can be viewed as an \s via its disjoint cycle decomposition; 
for definiteness, we define a \emph{cycle} to be a list whose smallest entry occurs first.

To \emph{delete} a run from a standard \s means to remove it and 
standardize what's left (replace smallest entry by 1, second smallest by 2, and so 
on). Thus deleting the run 23 from 178/23/465/9 yields 156/243/7. To 
\emph{insert} a run $i+1,\,\ldots, i+j$ into a standard \s $\pi$ 
means to increment by $j$ all elements of $\pi$ that exceed $i$ and adjoin 
$i,i+1,\ldots, j$ as a new block. The result will be a standard \s 
provided  $0 \le i \le \v\, \pi \,\v$. For example, inserting the run 456 into 15/342 yields 
 18/372/456. When runs are successively 
deleted from an \s, the order of deletion is immaterial and 
the result is always the same run-free \s. For 178/23/465/9, the result is 156/243.

Let $\p$ denote the set of all standard \ss, 
including the empty one $\ep$.  A \emph{run-closed} family $\fcal$ of 
\ss is a subset of \p that is closed under insertion and deletion of runs.

Some examples of run-closed families and, where available, their counting sequences are

\vspace*{-4mm}

\begin{itemize}
    \item  the family $\p$ itself \cite{sets08}  \htmladdnormallink{A000262}{http://oeis.org/A000262}
 \vspace*{-2mm}
    \item  set partitions, \htmladdnormallink{A000110}{http://oeis.org/A000110}
 \vspace*{-2mm}
    \item  noncrossing \ss  \htmladdnormallink{A088368}{http://oeis.org/A088368}
 \vspace*{-2mm}
    \item  nonoverlapping \ss 
 \vspace*{-2mm}
    \item  permutations, via the disjoint cycle decomposition  \htmladdnormallink{A000142}{http://oeis.org/A000142}
  \vspace*{-2mm}  
    \item  the intersection of any collection of run-closed families.    
   \end{itemize}  

\vspace*{-4mm}
   
   An ordinary set partition is an \s in which each block is an 
increasing list, and it can be represented graphically as the numbers 
$1,2,\ldots,n$ arranged around a circle with a line joining each pair 
of entries that are in the same block. It is noncrossing if no two 
lines cross. The run-closed property of noncrossing partitions is 
evident from this representation. Similarly, a set partition is 
nonoverlapping if the lines joining the smallest and largest entry of 
each block are noncrossing, a property that is also preserved under 
insertion/deletion of runs. An \s is 
\emph{noncrossing} if its underlying partition is noncrossing. We will have 
more to say about the last example in Section \ref{canonical}.

On the other hand, the family of nonnesting partitions is not 
run-closed. A partition is nonnesting if there is no 
quadruple $a<b<c<d$ with $a,d$ both in one block and $b,c$ both in another. 
Inserting the run 23 into the one-block nonnesting 
partition 12 produces the nesting partition 14/23. 

Now we can state our result for \ss, proved in the next two sections.
\begin{theorem}\label{main}
    Let $\fcal$ be a run-closed family of \ss with size \gf $f(x)$. 
    Then the run transform $F(x,y)$ of $f(x)$ counts $\fcal$ by size 
    and number of runs.
\end{theorem}

\section{Run-closed families with a singleton basis}\label{basis}
\vspace*{-5mm}
A run-closed family \q of \ss is 
determined by its \emph{run-free} members. This is because all members of \q can be obtained by successively 
inserting runs into its run-free members. We call the set of 
run-free \ss in a run-closed family $\q$ the \emph{basis} of $\q$.

Every set of run-free \ss  in
\p serves as a basis for a run-closed family of \ss.
We have the following two easily proved results for ordinary 
partitions.
\begin{lemma} 
   A standard noncrossing partition is either empty or contains a run.  \qed
\end{lemma} 
\begin{cor} \label{f0NC}
  The singleton set 
  consisting of the empty partition is the basis for the family of 
  noncrossing partitions.
  \qed
\end{cor}

 Every \s  can be 
  successively pruned of runs from right to left, leaving a run-free \s 
  (possibly empty) and a sequence of runs, its \emph{run-deletion sequence},
  from which the original \s can be recovered, as illustrated.
  \[
  \begin{array}{ll}
      \textrm{current \s} & \textrm{deleted run}  \\ 
      \hline
      1\ 12\ 10\ / \ 2\ 6\ 8\ / \ 3\ / \ 4\ 5\ / \ 7\ / \ 9\ 11\quad & 7  \\
      1\ 11\ 9\ / \ 2\ 6\ 7\ / \ 3\ / \ 4\ 5\ / \ 8\ 10 & 4\ 5  \\
      1\ 9\ 7\ / \ 2\ 4\ 5\ / \ 3\ / \ 6\ 8 &  3  \\
      1\ 8\ 6\ / \ 2\ 3\ 4\ / \ 5\ 7 & 2\ 3\ 4 \\
      1\ 5\ 3\ / \ 2\ 4 & \\ \hline
  \end{array}  
  \]
\centerline{ run-free \s =  1\,5\,3\,/\,2\,4,\quad run-deletion 
sequence = ( 2\,3\,4, 3, 4\,5, 7) }    
 
Furthermore, the number of runs in the original \s is 
captured in the run-deletion sequence as the number of runs that are 
disjoint from their immediate predecessor. (The first run vacuously 
meets this condition.)
This is because, in reconstructing the \s, when a new run is 
inserted, the only existing run that it can destroy is its 
predecessor run (if present) and it will do so precisely when 
it overlaps its predecessor.  A run-deletion sequence is of course 
specified by the first entries and lengths of its members, say 
$(a_{i})_{i=1}^{r}$ and $(\ell_{i})_{i=1}^{r}$ in reverse order of deletion. In the example 
$(a_{i})_{i=1}^{4}=(2,3,4,7);\ (\ell_{i})_{i=1}^{4}=(3,1,2,1)$.

For a run-closed family $\fcal$ of \ss, let 
$\fcal(n) =\{\rho\in \fcal\,:\vert\, \rho\,\vert =n\}$,
the members of $\fcal$ of size $n$.

\begin{prop} \label{charRunDeletion} Fix a run-free \s $\pi$ of size $k$.
   Let $\fcal$ denote the set of \ss that prune to $\pi$, and 
   suppose $n>k$.
Then the run-deletion sequences of \ss in  $\fcal(n)$, as 
specified by $(a_{i})_{i=1}^{r}$ and $(\ell_{i})_{i=1}^{r}$, 
are characterized by the following conditions:
\[
  \begin{array}{l}
      \textrm{\emph{$r\ge 1$ and all $a$'s and $\ell$'s are positive integers}},  \\
      k+\ell_{1}+\ell_{2}+\ldots +\ell_{r}=n, \\
      a_{1}<a_{2}< \ldots < a_{r},\\
      a_{1} \le k+1,  \\
      a_{2} \le k+\ell_{1}+1,  \\
      a_{3} \le k+\ell_{1}+\ell_{2}+1,  \\
       \quad   \vdots  \\
      a_{r} \le k+\ell_{1}+\ell_{2}+\ldots +\ell_{r-1}+1. 
  \end{array}
\]
\end{prop} 

Proof. The first two conditions are obvious. 
Now, when a run is deleted, the 
result is still a standard \s. Clearly, $a_{r}+\ell_{r}\le n+1$ and so $a_{r} \le 
n-\ell_{r}+1 = k+\ell_{1}+\ell_{2}+\ldots +\ell_{r-1}+1$, and similarly 
for the other inequalities.
Because runs are deleted right to left, we get  $a_{1}<a_{2}< \ldots < 
a_{r}$. 

Conversely, when runs are inserted successively into the run-free \s $\pi$ to build 
up members of $\fcal(n)$, the runs are arbitrary subject only to the 
conditions that the run currently being inserted begins at an integer 
no larger than 1 + the size of the \s it's being inserted into, 
for otherwise there would be a gap and the resulting \s would 
not have an initial segment of the positive integers as support. \qed

Now we can establish
\begin{prop}\label{onlysize}
Fix a run-free standard \s $\pi$. The number of \ss of given size and run count that 
prune to $\pi$ depends only on the size of $\pi$, not on its actual 
blocks.    
\end{prop}
Proof. \quad 
Suppose $\g_{1}$ is a run-closed family all of whose 
members prune to $\pi_{1}$ and $\g_{2}$ is a run-closed family all of whose 
members prune to $\pi_{2}$. Suppose further that $\pi_{1}$ and 
$\pi_{2}$ have the same size $k$. We wish to show that $\v\, 
\g_{1}(n)\,\v = \v\, \g_{2}(n)\,\v$ for all $n>k$ (it's obviously true 
for $n=k$). Since the characterization of the run-deletion sequences for $\g_{1}(n)$ makes no 
reference to $\pi_{1}$ other than through its size, this equality follows from two observations:

\vspace*{-2mm}

(i) members of $\g_{i}(n)$ are in 1-to-1 correspondence with their 
run-deletion sequences, $i=1,2$.

\vspace*{-2mm}

(ii) the set of run-deletion sequences for $\g_{1}(n)$ is 
\emph{identical} with that for $\g_{2}(n)$: both sets are characterized 
by the conditions of Proposition \ref{charRunDeletion}. \qed

\section{ A canonical family with singleton bases}\label{canonical}
\vspace*{-5mm}
Here we consider a canonical singleton basis of each size $k$, say the one-block \s 
$k,k-1,\ldots,2,1$ (the empty \s if $k=0$), and let $\fcal_{k}$ denote the run-closed 
family generated by it. A technical difficulty arises if $k=1$ because the only \s of 
size 1 is the one-block \s, (1), which is a run. 
(A workaround would be to color it red, maintain the color when runs are inserted---inserting 
12 into \red{1} yields 12/\red{3}---and declare that red entries are not to be considered runs.)

By Corollary \ref{f0NC}, $\fcal_{0}$ is the family of noncrossing 
ordinary partitions, and
$\fcal_{k}$ can be characterized for general $k$: it consists of the  noncrossing \ss 
in which (i) all blocks are increasing except for one block of length 
$k$ which is decreasing and (ii) this decreasing block is not covered by entries in another block, 
in other words, there is no increasing block containing integers $a<b$ such that 
all entries of the decreasing block lie in the interval $[a,b]$.   

For $k\ge 2$, a slight modification of the Simion 
bijection gives a bijection from nonnegative paths that start at 
$(0,k)$ to $\fcal_{k}$. First, prepend $k$ upsteps to turn the nonnegative path 
into a Dyck path.

Apply the Simion map of Section \ref{catconv}, as illustrated on an example with $k=3.$
The blocks are in the natural order and each block is decreasing.
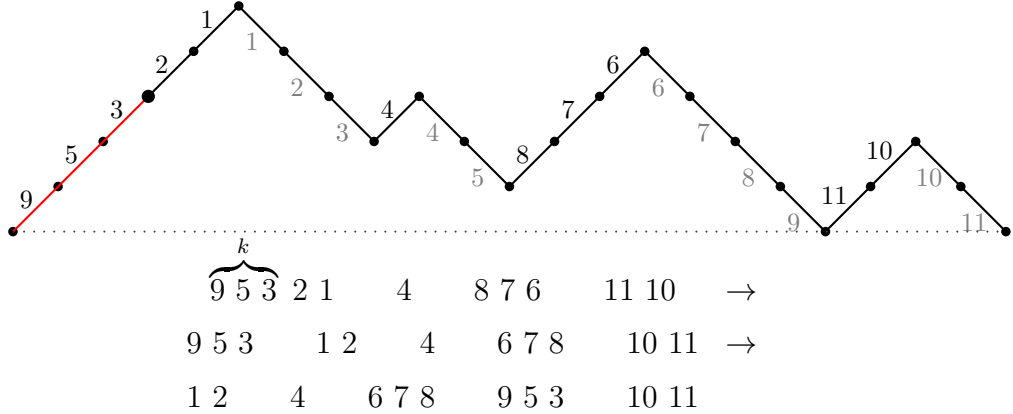
\begin{center} 
\begin{pspicture}(-8.5,0.5)(13,3)
\psset{unit=.6cm}     

\psdots(-11,0)(-10,1)(-9,2)(-8,3)(-7,4)(-6,5)(-5,4)(-4,3)(-3,2)(-2,3)(-1,2)
\psdots(0,1)(1,2)(2,3)(3,4)(4,3)(5,2)(6,1)(7,0)(8,1)(9,2)(10,1)(11,0)

\psline[linecolor=red](-11,0)(-8,3)
\psline(-8,3)(-6,5)(-3,2)(-2,3)(0,1)(3,4)(7,0)(9,2)(11,0)

\rput(-10.7,0.7){\textrm{{\footnotesize 9}}}
\rput(-9.7,1.7){\textrm{{\footnotesize 5}}}
\rput(-8.7,2.7){\textrm{{\footnotesize 3}}}
\rput(-7.7,3.7){\textrm{{\footnotesize 2}}}
\rput(-6.7,4.7){\textrm{{\footnotesize 1}}}

\rput(-2.7,2.7){\textrm{{\footnotesize 4}}}
\rput(0.3,1.7){\textrm{{\footnotesize 8}}}
\rput(1.3,2.7){\textrm{{\footnotesize 7}}}
\rput(2.3,3.7){\textrm{{\footnotesize 6}}}

\rput(7.2,0.8){\textrm{{\footnotesize 11}}}
\rput(8.2,1.8){\textrm{{\footnotesize 10}}}

\gray1{\rput(-3.7,2.2){\textrm{{\footnotesize 3}}}  
\rput(-4.7,3.2){\textrm{{\footnotesize 2}}} 
\rput(-5.7,4.2){\textrm{{\footnotesize 1}}}
\rput(-0.7,1.2){\textrm{{\footnotesize 5}}}
\rput(-1.7,2.2){\textrm{{\footnotesize 4}}}
\rput(3.3,3.2){\textrm{{\footnotesize 6}}}
\rput(4.3,2.2){\textrm{{\footnotesize 7}}}
\rput(5.3,1.2){\textrm{{\footnotesize 8}}}
\rput(6.3,0.2){\textrm{{\footnotesize 9}}}
\rput(9.3,1.2){\textrm{{\footnotesize 10}}}
\rput(10.3,0.2){\textrm{{\footnotesize 11}}}}

\newgray{darkgray}{.25}
\psline[linecolor=darkgray,linestyle=dotted](-11,0)(11,0)

\psset{dotsize=5pt 0}
\psdots(-8,3)
\end{pspicture}
\end{center}
\[
\begin{array}{cc}
    \overbrace{9\ 5\ 3}^{{k}} \ 2 \ 1\qquad 4  \qquad  8\ 7\ 6 \qquad  
    11\ 10 & \rightarrow \\[1mm]
    9\ 5\ 3 \qquad  1 \ 2   \qquad  4  \qquad  6\ 7\ 8  \qquad   10\ 11 
     & \rightarrow \\[1mm]
     1 \ 2   \qquad  4  \qquad  6\ 7\ 8  \qquad  9\ 5\ 3 \qquad   10\ 11
\end{array}
\]

\noindent The first block will have length $\ge k$ and will end with 1, and all 
blocks will be decreasing. Split the first block after the $k$-th 
entry into 2 blocks. Reverse all blocks except the new first block, 
and transfer the new first block to the appropriate position so that 
first entries are increasing. The resulting \s is a member 
of $\fcal_{k}$ and the mapping is reversible.

This correspondence preserves size and identifies runs with pyramid 
ascents. So we have
\begin{prop}\label{Fk}
    The \gf for $\fcal_{k}$ by size and number of runs is 
    $x^{k}C(x,y)^{k+1}$.
\end{prop}

The linearity of the run transform, along with Propositions \ref{onlysize} and \ref{Fk},
now yields 
\begin{prop} \label{basic}
    Let $\fcal$ be any run-closed family of \ss  
    containing $a_{k}\ (\ge 0)$ run-free \ss of size 
    $k,\ k\ge 0$.
    Then $F(x,y):=\sum_{k\ge 0}a_{k}x^{k}C(x,y)^{k+1}$ is the \gf to count $\fcal$ by 
    size and number of runs. \qed
\end{prop}
  
\noindent \textbf{Proof of Theorem \ref{main}}.\quad  
In the notation of Proposition \ref{basic}, 
the \gf by size of the run-closed family $\fcal$ is 
\[
f(x)=F(x,1)=\sum_{k\ge 0}a_{k}x^{k}C(x)^{k+1},
\]
since $C(x,1)=C(x)$. Lemma \ref{prop1} then yields that the 
run transform of $f(x)$ is indeed $F(x,y)$.

\textbf{Remark}. Theorem 8 applied to the case of set partitions implies that the
bivariate generating function for set partitions according to size and number of runs is
$\sum_{n\ge 0} B_n x^n (1-x)^{n+1}/(1 - x y)^{n+1}$, where $B_n$ are the Bell numbers. 
In particular,
the generating function for the number of partitions such that no block is a run is
$(1-x)\sum_{n\ge 0} B_n \big(x(1-x)\big)^{n}$ \cite[Exercise 111,\,pp. 137,\,192--3]{ec1}.

\section[Generalization of Theorem 8]{Generalization of Theorem \ref{main}}\vspace*{-5mm}
Fix a positive integer $j$. A $j$-compatible \s 
is one in which each block has length divisible by $j$. Define its 
size to be $n/j$ where $n$ is the cardinality of its support set. A 
$j$-compatible run ($j$-run for short) is one whose length and last 
entry are both divisible by $j$. A $j$-run-closed family of 
$j$-compatible \ss is one that is closed under insertion and 
deletion of $j$-runs.
\begin{theorem}
    \label{mainGeneral}
    Let $\fcal$ be a $j$-run-closed family of $j$-compatible \ss with size 
    \gf $f(x)$. 
    Then the run transform $F(x,y)$ of $f(x)$ counts $\fcal$ by size 
    and number of $j$-runs.
\end{theorem}
Proof. 
The ``$j$'' analogue of $\fcal_{k}$ is the family of $j$-compatible \ss 
with singleton basis $(jk,jk-1,\ldots,1)$, which corresponds under Simion's 
bijection to the family of $(j,jk,0)$-$U$-$D$ paths.
This bijection preserves size and identifies $j$-runs with $j$-nice 
pyramid ascents. Apply Theorem \ref{jmdUD}. \qed

As an example, we have the following result.

\begin{cor}
If $f(x)$ denotes the \gf for permutations of $[2n]$ in which each cycle
has even length 
(\htmladdnormallink{A001818}{http://oeis.org/A001818})
by size $n$, then the run transform of $f$ counts these permutations by
size and by number of cycles that consist of consecutive integers 
ending at an even integer.
\end{cor}


\section[U-F-D paths]{$\mbf{U}$-$\mbf{F}$-$\mbf{D}$ paths}\vspace*{-5mm}

Fix nonnegative integers $m$ and $d$ and consider 
the class $\a_{m,d}$ of lattice paths of upsteps $U=(1,1)$, downsteps 
$D=(1,-1)$, and flatsteps $F=(2,0)$ that 
start at $(0,m)$, end on the
$x$-axis and that reach lowest level $-d$, with size measured by (number of 
flatsteps) + (number of downsteps). Thus $\a_{0,0}$ is the 
class of Schr\"{o}der paths with the usual measure of size. The 
definition of pyramid ascent carries over to $\a_{m,d}$. 
\begin{lemma}\label{lem1}
    Let $F(x,y,z_{0},z_{1},z_{2},...)$ denote the \gf for Schr\"{o}der paths with $x$ marking size, 
    $y$ marking number of pyramid ascents, and $z_{i}$ marking  number of flatsteps 
    at level $i,\ i\ge 0$.
    Thus $f(x,z$'s)$:= F(x,1,z$'s) is the \gf disregarding pyramid ascents. 
    Then the run transform of $f$ is $F$.
\end{lemma}
Proof.
Set $F_{j} = F(x,y,z_{j},z_{j+1},z_{j+2},...)$. Thus $F_{0}=F$.
A Schr\"{o}der path is either empty or, by considering the first non-upstep, 
begins with one of the 
prefixes $F,\ UD,$ $U^{r}D\ (r\ge 2),\ U^{r}F\ (r\ge 1)$. 
Thus a nonempty Schr\"{o}der path has precisely one 
of the following forms, illustrated for $r=3$, where the $S_{i}\ (i\ge 
0)$ denote arbitrary Schr\"{o}der paths.

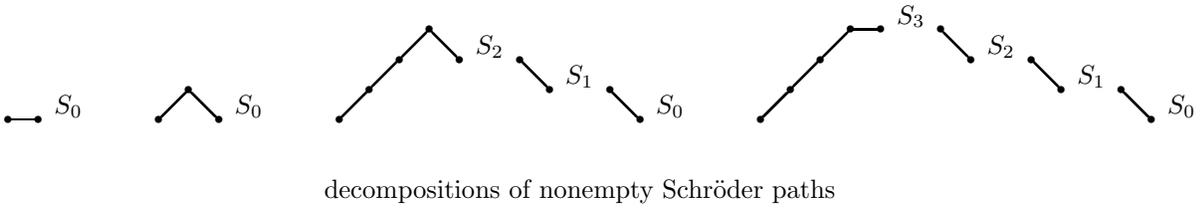
\Einheit=0.4cm
\[
\Pfad(-20,0),1\endPfad
\Pfad(-15,0),34\endPfad
\Pfad(-9,0),3334\endPfad
\Pfad(-3,2),4\endPfad
\Pfad(0,1),4\endPfad
\Pfad(5,0),3331\endPfad 
\Pfad(11,3),4\endPfad
\Pfad(14,2),4\endPfad
\Pfad(17,1),4\endPfad
\Label\o{{\textrm{\footnotesize $S_{0}$}}}(-18,-0.2)
\Label\o{{\textrm{\footnotesize $S_{0}$}}}(-12,-0.2)
\Label\o{{\textrm{\footnotesize $S_{2}$}}}(-4,1.8)
\Label\o{{\textrm{\footnotesize $S_{1}$}}}(-1,0.8)
\Label\o{{\textrm{\footnotesize $S_{0}$}}}(2,-0.2)
\Label\o{{\textrm{\footnotesize $S_{3}$}}}(10,2.8)
\Label\o{{\textrm{\footnotesize $S_{2}$}}}(13,1.8)
\Label\o{{\textrm{\footnotesize $S_{1}$}}}(16,0.8)
\Label\o{{\textrm{\footnotesize $S_{0}$}}}(19,-0.2)
\DuennPunkt(-20,0)
\DuennPunkt(-19,0)
\DuennPunkt(-15,0)
\DuennPunkt(-14,1)
\DuennPunkt(-13,0)
\DuennPunkt(-9,0)
\DuennPunkt(-8,1)
\DuennPunkt(-7,2)
\DuennPunkt(-6,3)
\DuennPunkt(-5,2)
\DuennPunkt(-3,2)
\DuennPunkt(-2,1)
\DuennPunkt(0,1)
\DuennPunkt(1,0)
\DuennPunkt(5,0)
\DuennPunkt(6,1)
\DuennPunkt(7,2)
\DuennPunkt(8,3)
\DuennPunkt(9,3)
\DuennPunkt(11,3)
\DuennPunkt(12,2)
\DuennPunkt(14,2)
\DuennPunkt(15,1)
\DuennPunkt(17,1)
\DuennPunkt(18,0)
\Label\o{\textrm{\footnotesize decompositions of nonempty Schr\"{o}der 
paths}}(-1,-3)
\]

\vspace*{2mm}
From this schematic picture, we see that
\begin{eqnarray*}
  F_{0}   & = &   1 + xzF_{0}  + xyF_{0} + \Big(\sum_{r\ge 2}(F_{r-1}\ldots F_{1}-1)F_{0} +\sum_{r\ge 
2}x^{r}yF_{0}\Big) + \sum_{r\ge 1}x^{r+1}z_{r}F_{r}F_{r-1}\ldots F_{0} 
\end{eqnarray*}
from which we obtain by routine manipulation
\begin{equation}
    \frac{1-xy}{1-x}F_{0} = 1+\sum_{r\ge 1}x^{r} (1+z_{r-1})F_{r-1}F_{r-2}\ldots F_{0}
    \label{Fxyz}\, ,
\end{equation}
a recursion for $F=F_{0}$ (bear in mind that $F_{1},F_{2},\ldots$ are 
merely abbreviations for functions derived from $F$). This recursion 
has a unique solution for $F$ because it determines the constant term 
and then the coefficients of $x,x^{2},\ldots$ in turn.

Set $f_{j}(x,z_{j},z_{j+1},\ldots) 
=F_{j}(x,1,z_{j},z_{j+1},\ldots)$. Thus $f_{0}=f.$ From (\ref{Fxyz}) with $y=1$, we have
\begin{equation}
f_{0}=1+\sum_{r\ge 1}x^{r} (1+z_{r-1})f_{r-1}f_{r-2}\ldots f_{0}.
 \label{fxz}
\end{equation}
The run transform of $f_{0}(x,z_{0},z_{1},\ldots)$ is 
\[
H_{0}(x,y,z_{0},z_{1},\ldots):=\frac{1-x}{1-xy}\,f_{0}\left(\frac{x(1-x)}{1-xy},z_{0},z_{1},\ldots\right),
\]
and we define $H_{j},\ j\ge 1 $ by relabeling $z$ indices just as for 
$F_{j}$.
To verify that $H_{0}$ and $F_{0}$ are equal, replace $x$ by 
$\frac{x(1-x)}{1-xy}$ in (\ref{fxz}) to obtain
\begin{multline}\notag
\frac{1-xy}{1-x}\left(\frac{1-x}{1-xy}\,f_{0}\Big(\frac{x(1-x)}{1-xy},z_{0},z_{1},\ldots\Big)\right) = \\
1+\sum_{r\ge 1}x^{r} (1+z_{r-1}) 
\frac{1-x}{1-xy}f_{r-1}\Big(\frac{x(1-x)}{1-xy},z\textrm{'}s\Big) \ldots \frac{1-x}{1-xy}f_{0}\Big(\frac{x(1-x)}{1-xy},z\textrm{'s}\Big)
\end{multline}
or
\begin{equation}
\frac{1-xy}{1-x}H_{0}=1+\sum_{r\ge 1}x^{r} (1+z_{r-1}) 
H_{r-1}H_{r-2}\ldots H_{0} 
\label{hxyz}
\end{equation}
Comparing (\ref{Fxyz}) and (\ref{hxyz})  we see that $H_{0}=F_{0}$  
because, as noted above, (\ref{Fxyz}) has a unique solution. \qed

\begin{theorem}
  Let $G(x,y,z_{-d},z_{-(d-1)},...,z_0,z_1,z_2,...)$ denote 
  the \gf for $\a_{m,d}$ with $x$ marking size, 
    $y$ marking  number of pyramid ascents, and $z_{i}$ marking  number of flatsteps 
    at level $i$. Thus $g(x,z$'s)$:= G(x,1,z$'s) is the \gf disregarding pyramid ascents. 
    Then the run transform of $g$ is $G$.
\end{theorem}
Proof. 
A path in  $\a_{j,d}$ has the form below,
illustrated for $d=2$, where $S_{i}$ and $S_{i}'$ denote Schr\"{o}der paths.
\vspace*{4mm}
\Einheit=0.5cm
\[
\hspace*{-95mm}
\red{\Pfad(-3,0),1111111111111111111111111\endPfad}
\red{\Pfad(-3,0),22222\endPfad}
\red{\Label\o{{\textrm{\footnotesize O}}}(-3,-1.2)
\DuennPunkt(-3,4)  }
\Pfad(-1,4),4\endPfad
\Pfad(3,3),4\endPfad
\Pfad(10,0),4\endPfad
\Pfad(13,-1),4\endPfad
\Pfad(7,1),4\endPfad
\Pfad(16,-2),3\endPfad
\Pfad(19,-1),3\endPfad
\SPfad(5,2),4\endSPfad
\Label\o{{\textrm{\footnotesize $S_{j}$}}}(-2,3.8)
\Label\o{{\textrm{\footnotesize $S_{j-1}$}}}(1.5,2.8)
\Label\o{{\textrm{\footnotesize $S_{0}$}}}(9,-0.1)
%
\Label\o{{\textrm{\footnotesize $S_{-1}$}}}(12,-1.2)
\Label\o{{\textrm{\footnotesize $S_{-2}$}}}(15,-2.2)
\Label\o{{\textrm{\footnotesize $S_{-1}'$}}}(18,-1.2)
\Label\o{{\textrm{\footnotesize $S_{0}'$}}}(21,-0.1)
\DuennPunkt(-1,4)
\DuennPunkt(0,3)
\DuennPunkt(3,3)
\DuennPunkt(4,2)
\DuennPunkt(7,1)
\DuennPunkt(8,0)
\DuennPunkt(10,0)
\DuennPunkt(11,-1)
\DuennPunkt(13,-1)
\DuennPunkt(14,-2)
\DuennPunkt(16,-2)
\DuennPunkt(17,-1)
\DuennPunkt(19,-1)
\DuennPunkt(20,0)
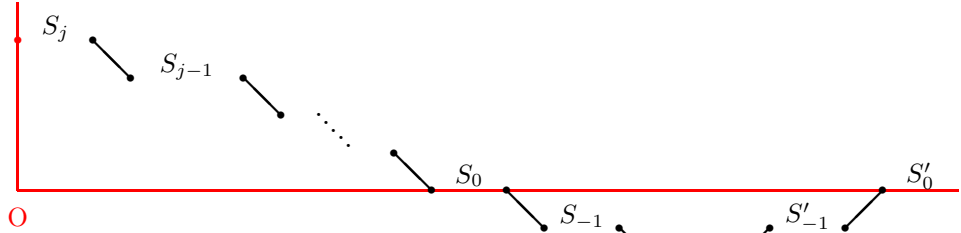
\Label\o{{\textrm{\footnotesize decomposition of path in 
$\a_{j,d}$}}}(9,-4.5)
\]
\vspace*{1mm}

Consequently,
\begin{equation}
G=x^{j}F_{j}F_{j-1}\ldots F_{0}  + \sum_{k = 1}^{d}x^{j+k} 
F_{j}F_{j-1}\ldots F_{0}F_{-1}\ldots F_{-k}
\prod_{i=0}^{k-1}\Big( \frac{H_{-i}}{y}+(F_{-i}-H_{-i}) \Big),
\label{Gxyz}
\end{equation}
where $H_{0}(x,y,z_{0},z_{1},\ldots)$ is the \gf for nonempty 
Schr\"{o}der  paths that start with a pyramid, and $H_{-i}\ (i\ge 1)$ is 
obtained from $H_{0}$ by relabeling $z$ indices in the usual way. The 
introduction of $H_{0}$ is necessary because if a Schr\"{o}der path 
$S_{-i}'$ ($i\ge 0$) begins with a pyramid, then the immediately preceding upstep 
kills the initial pyramid ascent in $S_{-i}'$.
Clearly,
\begin{equation}
H_{0}(x,y,z_{0},z_{1},\ldots)=\sum_{k\ge 1}x^{k} y F = \frac{xyF}{1-x}.
\label{Hxyz}
\end{equation}
Now let $g(x,z_{-d},z_{-(d-1)},\ldots, z_{0},\ldots) 
=G(x,1,z_{-d},z_{-(d-1)},\ldots, z_{0},\ldots)$.
Thus
\begin{equation}
g = x^{j}f_{j}f_{j-1}\ldots f_{0}  + \sum_{k = 1}^{d}x^{j+k} 
f_{j}f_{j-1}\ldots f_{0}f_{-1}\ldots f_{-k}
\prod_{i=0}^{k-1}f_{-i}.
\label{gxz}
\end{equation}

From (\ref{Hxyz}), we have
\begin{equation}
\frac{H_{-i}}{y}+(F_{-i}-H_{-i}) = \frac{1-xy}{1-x}F_{-i}.
\label{hf}
\end{equation}

Using (\ref{gxz}), the run transform of $g$ is
\begin{multline}\notag
\frac{1-x}{1-xy}g\Big(\frac{x(1-x)}{1-xy},z\textrm{'s}\Big) = \\[1mm]
x^{j} \frac{(1-x)^{j+1}}{(1-xy)^{j+1}}f_{j}\Big(\frac{x(1-x)}{1-xy},z\textrm{'s}\Big) \ldots
f_{0}\Big(\frac{x(1-x)}{1-xy},z\textrm{'s}\Big) + \\[1mm]
\sum_{k=1}^{d}x^{j+k} \frac{(1-x)^{j+k+1}}{(1-xy)^{j+k+1}} 
f_{j}\Big(\frac{x(1-x)}{1-xy},z\textrm{'s}  \Big) f_{j-1}\Big(\frac{x(1-x)}{1-xy},z\textrm{'s}  \Big)
\ldots f_{0}\Big(\frac{x(1-x)}{1-xy},z\textrm{'s}  \Big) \times 
\\[1mm]
f_{-1}\Big(\frac{x(1-x)}{1-xy},z\textrm{'s}  \Big) \ldots f_{-k}\Big(\frac{x(1-x)}{1-xy},z\textrm{'s}  \Big) 
\prod_{i=0}^{k-1}f_{i}\Big(\frac{x(1-x)}{1-xy},z\textrm{'s}  \Big)
\end{multline}

\vspace*{-10mm}

\[
= x^{j} F_{j}\ldots F_{0} + \sum_{k=1}^{d}x^{j+k}F_{j}F_{j-1}\ldots 
F_{0}F_{-1}\ldots F_{-k} \prod_{i=0}^{k-1}\frac{1-xy}{1-x}F_{-i},
\]
which, in view of (\ref{hf}), is the same expression as in (\ref{Gxyz}). 
The run transform of $g$ is thus $G$. \qed

\section{Concluding remark}\vspace*{-5mm}
We believe there is a version of our results that includes both $j$-upsteps
and flatsteps. The setting 
is paths of $j$-upsteps $U_{j}=(j,j)$, flatsteps $F=(2,0)$, and downsteps 
$D=(1,-1)$, that start at $(0,m)$ and end on the $x$-axis; $j$ and $m$ 
nonnegative integers. In this generality, the size of a path is
$\lfloor (number of D$'s$ + number of F$'s$)/j \rfloor$ (so it doesn't matter 
whether we consider flatsteps to be of length 1 or
2). Furthermore, the ``nice'' pyramid ascents to count are 
those whose endpoint $(a,b)$ has $a$
divisible by $j$ rather than those $b$ divisible by $j$; that is, the abscissa rather than the
ordinate of the endpoint is divisible by $j$. These conditions are equivalent 
if $j=1$ or if there are no flatsteps and $m$ is divisible by $j$.
\begin{conj} Fix nonnegative integers $j$ and $m$.
  Let $G(x,y,z_{-d},z_{-(d-1)},...,z_0,z_1,z_2,...)$ denote 
  the \gf for $U_{j}$-$F$-$D$ paths with $x$ marking size, 
    $y$ marking  number of ``nice'' pyramid ascents, and $z_{i}$ marking  number of flatsteps 
    at level $i$. Thus $g(x,z$'s)$ \ := G(x,1,z$'s) is the \gf disregarding pyramid ascents. 
    Then the run transform of $g$ is $G$.
\end{conj}

\end{document}